\documentclass[opre,nonblindrev]{informs3} 

\OneAndAHalfSpacedXI 

\usepackage{amssymb, amsmath, amsfonts}
\usepackage{mathrsfs, color, bm, enumitem, url, xr}
\usepackage{graphicx, transparent} 
\usepackage{multirow, titlecaps, accents, epstopdf}
\usepackage{caption}
\usepackage{subcaption}
\usepackage{alphalph,etoolbox}
\patchcmd{\subequations}{\alph{equation}}{\alphalph{\value{equation}}}{}{}

\graphicspath{{figures/}}

\usepackage{algorithmic}            
\usepackage{algorithm}

\usepackage[utf8]{inputenc}    
\usepackage{hyperref}          
\usepackage{url}               
\usepackage{booktabs}          
\usepackage{amsfonts}          
\usepackage{nicefrac}          
\usepackage{microtype}         



\makeatletter
\def\munderbar#1{\underline{\sbox\tw@{$#1$}\dp\tw@\z@\box\tw@}}
\makeatother

\newcommand{\bea}{\begin{equation*}\begin{aligned}}
\newcommand{\eea}{\end{aligned}\end{equation*}}

\DeclareMathOperator{\st}{s.t.}

\newcommand{\set}[1]{\mathcal{#1}}
\newcommand{\R}{\mathbb{R}}

\usepackage{url}
\usepackage{natbib}
\usepackage{comment} 
\usepackage{xcolor}
\usepackage{anyfontsize}

\usepackage{appendix}

\bibpunct[, ]{(}{)}{,}{a}{,}{,}%
\def\bibfont{\small}%
\setcitestyle{citesep={;}}

\DeclareMathOperator*{\dom}{dom}
\TheoremsNumberedThrough     
\ECRepeatTheorems

\EquationsNumberedThrough    

\usepackage{xcolor}

\usepackage[normalem]{ulem}

\hypersetup{
	colorlinks=true,
	citecolor=blue,
	urlcolor=blue,
	linkcolor=blue,
	pdfborder={0 0 0},
}

\begin{document}
	
	
\RUNTITLE{Global Nonconvex Optimization with Integer Variables}
\TITLE{Global Nonconvex Optimization with Integer Variables}

\ARTICLEAUTHORS{
 		\AUTHOR{Dimitris Bertsimas}
		\AFF{Sloan School of Management, Massachusetts Institute of Technology, United States,
  \EMAIL{dbertsim@mit.edu}}
		\AUTHOR{Danique de Moor}
		\AFF{Operations Research Center, Massachusetts Institute of Technology, United States, 
			\EMAIL{demoor@mit.edu}}
\AUTHOR{Thodoris Koukouvinos}
		\AFF{Operations Research Center, Massachusetts Institute of Technolgy, United States, 
			\EMAIL{tkoukouv@mit.edu}}
 	
 		\AUTHOR{Demetrios Kriezis}
		\AFF{Massachusetts Institute of Technology, United States, 
			\EMAIL{dkriezis@mit.edu}}
  }

	\ABSTRACT{Nonconvex optimization refers to the process of solving problems whose objective or constraints are nonconvex. Historically, this type of problems have been very difficult to solve to global optimality, with traditional solvers often relying on approximate solutions. \cite{bertsimas2023rptbb} introduce a novel approach for solving continuous nonconvex optimization problems to provable optimality, called the Relaxation Perspectification Technique - Branch and Bound (RPT-BB). In this paper, we extend the RPT-BB approach to the binary, mixed-binary, integer, and mixed-integer variable domains. We outline a novel branch-and-bound algorithm that makes use of the Relaxation Perspectification Technique (RPT), as well as binary, integer, and eigenvector cuts. We demonstrate the performance of this approach on four representative nonconvex problems, as well as one real-world nonconvex optimization problem, and we benchmark its performance on BARON and SCIP, two state-of-the-art optimization solvers for nonconvex mixed-integer problems.
    Our results show that our method stands well against BARON, and often outperforms BARON, in terms of computational time and optimal objective value. Moreover, our results show that, while SCIP continues to lead as the state-of-the-art solver, the proposed algorithm demonstrates strong performance on challenging instances, successfully solving problems to global optimality that SCIP and BARON are unable to solve within the time limit.
	} 
	
	\KEYWORDS{}
	
	\maketitle
    

\section{Introduction}
Many real-world decision-making problems give rise to optimization models that are nonconvex for multiple reasons: they often involve integer variables as well as nonconvex objective functions and/or constraints defined over both continuous and discrete variables. Integer variables are essential when modeling discrete choices, such as selecting which facilities to open, which investments to make, or assigning tasks to resources, and are commonly used in applications including but not limited to facility location, healthcare operations, telecommunications network design, transportation and logistics, and vehicle routing. Meanwhile, nonconvex functions and/or constraints often arise from practical considerations such as economies of scale (producing more leads to lower average cost per unit) \citep{alkaabneh2019lagrangian}, concave revenue functions (selling more leads to diminishing marginal returns in revenue) \citep{huang2024marketing}, or startup and ramping costs \citep{byers2023long}, which can affect both continuous and discrete aspects of the model. These sources of nonconvexity frequently arise in applications such as energy systems \citep{byers2023long}, transportation \citep{alkaabneh2019lagrangian}, supply chain operations \citep{alkaabneh2019lagrangian}, and portfolio optimization under nonlinear risk measures \citep{ahmadi2019portfolio}.

Traditionally, optimization solvers have been highly effective for convex problems but face significant challenges when tackling nonconvex objectives. One reason is that nonconvex problems can have many local optima, making it hard for solvers to find the best possible solution. The presence of integer variables eliminates the possibility of relying on gradient-based local search methods, often leveraged in continuous nonconvex solvers. In addition, convex relaxations become less effective because the fractional solutions they produce do not satisfy the integer constraints of the original problem.

State-of-the-art global optimization solvers such as BARON \citep{baron} and SCIP \citep{scip} are capable of tackling a broad class of nonconvex optimization problems with integer variables, using sophisticated methods like branch-and-bound and cutting planes combined with convex relaxations. \cite{bertsimas2023rptbb} introduce the Reformulation Perspectification Technique with Branch and Bound (RPT-BB), an approach originally designed for solving nonconvex continuous problems where the nonconvexity arises from sums of products of convex functions. \cite{Hartman59} shows that every twice differentiable function on a compact convex set $\set X$ can be written as an SLC function. Moreover, from the Weierstrass Approximation Theorem, every nonconvex continuous function on a compact convex set $\set X$, can be approximated by a twice differentiable function with arbitrarily small error $\varepsilon$. Hence, RPT-BB can be directly applied to general continuous nonconvex optimization problems. The RPT-BB approach embeds the convex relaxation technique RPT into a spatial branch and bound framework that leverages the eigenvector branching strategy proposed in \cite{anstreicher2022solving}. It has shown remarkable strength in producing tight convex relaxations, and as a  result can efficiently compute the global optimum.

In this paper, we extend the RPT-BB method to a broad class of mixed-integer nonconvex optimization problems. We first relax all integer variables to the continuous domain. We then apply RPT to construct a tight convex relaxation of the original problem. Finally, we propose a hybrid branch-and-bound strategy that combines eigenvector-based branching with traditional integer branching rules to effectively navigate both continuous and discrete decision spaces.

\medskip \noindent 
\textbf{Contributions}. Our main contributions can be summarized as follows:
 \begin{enumerate}
\item We propose a way to obtain a convex relaxation from a nonconvex mixed-integer optimization problem, leveraging the RPT method. In this way we can obtain strong lower bounds, which can lead to the global optimum, when incorporated in the RPT-BB method.
\item We demonstrate the effectiveness of our approach by conducting numerical experiments on three convex maximization problems, a quadratically constrained quadratic optimization problem, and a dike height optimization problem. We show that our approach stands well against the state-of-the-art solvers SCIP and BARON, and is able to find the global optimum in several instances for which BARON and SCIP cannot prove global optimality within the time limit. 
\item We investigate three branching strategies within our framework and provide insights into their relative performance.
\end{enumerate}

\vspace{0.3cm}
\noindent  This paper is structured as follows: In Section \ref{sec:generic_problem}, we describe the generic mixed-integer nonconvex optimization problem that we consider. In Section \ref{sec:method}, we describe our method to find the global optimal solution of the considered problem. In Section \ref{sec:numexp}, we report numerical results, comparing the performance of our approach against BARON and SCIP, and finally in Section \ref{sec:conclusions} we summarize our conclusions. 


\medskip\noindent
	{\bf Notation.}
We use bold lowercase letters (e.g. $\bm a \in \mathbb{R}^n$) to denote column vectors and bold uppercase letters (e.g. $\bm A \in \mathbb{R}^{m \times n}$) for matrices. The $i$-th entry of $\bm a$ is written $a_i$, the $i$-th column of $\bm A$ is denoted by $\bm A_i$, and the entry of $\bm A$ in the $i$-th row and $j$-th column is denoted by $A_{ij}$. Calligraphic letters and their corresponding Roman letters denote finite index sets and their cardinalities respectively (e.g. $\set I = \{1, \ldots, I\}$). 
 The {\it perspective function} of a proper, closed and convex function $f:\R^{n_\nu} \rightarrow (-\infty, +\infty]$ is defined as $h(\bm \nu,t) = tf(\bm \nu/t)$ if $t>0$, and $h(\bm \nu,0) = \delta^*_{\dom (f^*)}(\bm \nu)$ for all $\bm \nu \in \mathbb{R}^{n_\nu}$ and $t \in \mathbb{R}_+$, where $\delta^*_{\set V}$ denotes the support function of the set $\set V$ and $f^*$ denotes the conjugate of $f$ defined through $f^*({{\bm w}}) = \sup_{\bm \nu}\left\{{{\bm \nu}}^\top{{\bm w}} - f({{\bm \nu}})\right\}$. 
 For ease of exposition, we use $ tf(\bm \nu/t)$ to denote the perspective function $h(\bm \nu,t)$ for the rest of this paper. 

\section{Generic problem formulation} \label{sec:generic_problem}
We consider a generic mixed-integer nonconvex optimization problem of the following form:	
\begin{equation}
    \begin{array}{cll}
        \displaystyle \min_{\bm x} &\quad f_0(\bm x) \\
        {\rm s.t.} 	&\quad  f_k(\bm x) \le 0, & \quad \forall k \in \mathcal{K}, \\
    &\quad \bm x \in \mathcal{X},
    \end{array}
\label{eq:genericprobleminteger}
\end{equation} 
	where $f_k : \mathbb{R}^{n_x} \rightarrow [-\infty, \infty]$ is a sum of linear times convex (SLC) function for all $k \in \mathcal{K}_0$, that is, 
    \begin{equation}\label{eq:SLCfunction}
        f_k(\bm x) = c_{0k} (\bm x) +  \sum_{i \in \set I_k}  (b_i - \bm a_i^T \bm x) c_{ik}(\bm x), 
    \end{equation}	and $c_{0k}, c_{ik}: \mathbb{R}^{n_x}\rightarrow (-\infty, + \infty]$, are proper, closed and convex functions for every $i \in \mathcal{I}_k$, $k\in \mathcal{K}_0$. The set $\mathcal{X}$ is defined by:
	\[ \mathcal{X} = \left\{ \bm x \equiv (\bm x_B, \bm x_D, \bm x_C) \mid \bm A^\top \bm x \le \bm b, \ \bm P^\top \bm x = \bm s, \ \bm h (\bm x) \le \bm 0 \right\},
	\]
\noindent where the vector $\bm{x}$ is expressed as the concatenation of the binary vector $\bm{x}_B \in \{0,1\}^{|\set{B}|}$, the integer vector $\bm x_D \in \mathbb{Z}^{|\set{D}|}$ and the continuous vector $\bm{x}_C \in \mathbb{R}^{|\set{C}|}$, with $n_x = |\set{B}| + |\set{D}| + |\set{C}|$ and $\set{B} \cap \set{C} \cap \set{D} = \emptyset$. Additionally, we have $\bm A \in \mathbb{R}^{n_x \times m_1}$, $\bm P \in \mathbb{R}^{n_x \times m_2}$, $\bm b \in \mathbb{R}^{m_1}$, $\bm s \in \mathbb{R}^{m_2}$, $\bm h (\bm x) = [h_0(\bm x) \ h_1(\bm x)  \ \cdots \ h_J (\bm x)]^\top$, and $h_j: \mathbb{R}^{n_x} \rightarrow (-\infty, + \infty ]$ are proper, closed and convex for every $j \in \mathcal{J}_0$. 
We make the following assumptions.
 \begin{assumption} \label{assumption:X}
		The set $\mathcal{X}$ is nonempty and compact.
\end{assumption}
\begin{assumption} \label{assumption:SCC}
    If $r_{ik}$ and $c_{ik}$ are both nonlinear, then~$r_{ik}(\bm x)\ge 0$ and $c_{ik}(\bm x) \ge 0 $ for all $\bm x \in \mathcal{X}$, for every $i \in \mathcal{I}_k$ and $k \in \mathcal{K}_0$.
    If $r_{ik}$ is linear and $c_{ik}$ is nonlinear, then $r_{ik}(\bm x)\ge 0$ for all $\bm x \in \mathcal{X}$, for every $i \in \mathcal{I}_k$ and $k \in \mathcal{K}_0$. If both $r_{ik}$ and $c_{ik}$ are linear, then we do not impose any assumption on these functions. 
\end{assumption}
 Moreover, in this section we will give some examples of SLC functions that often occur in the literature and in practice.

\begin{example}[Mixed integer nonconvex quadratic problem]
We have the following problem
\begin{equation}
    \begin{array}{cll}
        \displaystyle \min_{\bm x} &\quad \bm x^\top \bm P \bm x + \bm q^\top \bm x + r \\
        {\rm s.t.} 	&\quad  \bm x \in \mathcal{X}, \\ 
        &\quad x_i \in \mathbb{R}, &\quad i \in \mathcal{C}, \\ 
        &\quad x_i \in \mathbb{Z}, &\quad i \in \mathcal{D},
    \end{array}
\label{eq:example_quadratic}
\end{equation} 
where $\mathcal{X}$ contains linear constraints and $\bm P$ is not necessarily positive semidefinite. One application of this problem is the Max-Cut problem \citep{billionnet2007using}, which appears in various domains such as binary classification tasks in machine learning and circuit partitioning in VLSI design. Let $G = (V,E)$ be a weighted undirected graph. The max-cut problem consists in finding a partition $\mathcal{S}, V \backslash \mathcal{S}$ of the set of vertices $V$ such that the total weight of the crossing edges is maximum. Let $w_{ij}$ denote the weight of edge $(i,j)$, and $x_i$ denote the decision variable that is equal to 1 if and only if vertex $i$ is in $\mathcal{S}$. The resulting problem formulation is 
\begin{equation}
    \begin{array}{cll}
        \displaystyle \min_{\bm x} &\quad - \sum_{i=1}^{V} \sum_{j \neq i} w_{ij} x_i x_j - \sum_{i=1}^V \left( \sum_{j \neq i} w_{ij} \right) x_i \\
        {\rm s.t.} 	&\quad  x_i \in \{0,1\}, &\quad i \in \set B.
    \end{array}
\label{eq:example_max_cut}
\end{equation} 
The mixed-integer nonconvex quadratic problem is treated in Section \ref{sec:qcqp}.
\hfill $\square$
\end{example} 

\begin{example}[Dike height optimization]
Dike height optimization is of crucial importance to the Netherlands since 60$\%$ of the Netherlands is flood prone. \cite{dike} develop a model that optimizes dike heightenings in the Netherlands. The optimal solution is shown to be periodic, that is, every $t$ years the dike is increased by the same amount. In practical settings, one might however want to deviate from this solution as it might be for example more practical to coordinate the heightening of multiple dikes simultaneously. \citet{bertsimas2023rptbb} consider a time truncated adaptation of the dike heightening problem, allowing for deviations from the periodic schedule of every $t$ years, i.e.,  
\small
\begin{align} 
            \begin{array}{cll}
            \displaystyle	\min_{\bm x, \bm h \in \mathbb{Z}_{+}^N } \underbrace{\sum_{k \in \mathcal{K}_0} \left(C + bx_k\right) \exp\left(\lambda h_k - \delta t_k\right)}_{\rm Investment \ costs} + \underbrace{ \sum\limits_{k \in \mathcal{K}_0} \frac{S_0}{\beta_{\delta}}\left(\exp\left(\beta_{\delta} t_{k+1}\right) - \exp\left(\beta_{\delta} t_k \right)\right)\exp\left(-\theta h_k \right)}_{\rm Expected \ damage \ costs} + \underbrace{\vphantom{\sum_{k \in K} f_k} \frac{S_0}{\delta}\exp\left(\beta_{\delta}T - \theta h_K \right)}_{\rm Future \ damage \ costs},
            \end{array}
\end{align}
\normalsize
where $x_k$ is the integer increment of the dike height at time $t_k$, $h_k$ is the increase in dike height after $t_k$ years, i.e., $h_k = \sum_{i=0}^k x_i$, $h_K = \sum_{k \in \mathcal{K}_0} x_k$ and $\beta_{\delta}, \delta, \theta, \lambda, b, C, T$ and $S_0$ are constants. The objective is an SLC function. We treat this problem in Section \ref{sec:dike}.
\hfill $\square$
\end{example} 

\section{Method}  \label{sec:method} 
In this section, we describe our approach to obtain the global optimal solution of Problem \eqref{eq:genericprobleminteger}. Our approach comprises three steps: First, we relax the integer variables to obtain a nonconvex continuous formulation of the original mixed-integer problem. Then, we apply RPT to the relaxed problem to construct a tight convex relaxation. Finally, we propose a hybrid branch and bound algorithm that combines the eigenvector-based branching of \cite{anstreicher2022solving} used in RPT-BB for continuous problems, with classical integer branching techniques used in mixed-integer problems. The steps are explained in more detail below.

\medskip\noindent
{\bf Step 1: Relaxation of integer variables.} We begin by relaxing the integrality constraints of the mixed-integer nonconvex optimization problem \eqref{eq:genericprobleminteger}, converting it into a continuous problem. The relaxed formulation is given by:
\begin{equation}
    \begin{array}{cll}
        \displaystyle \min_{\bm x} &\quad f_0(\bm x) \\
        {\rm s.t.} 	&\quad  f_k(\bm x) \le 0, & \quad \forall k \in \mathcal{K}, \\
    &\quad \bm x \in \set X_C, 
    \end{array}
\label{eq:genericproblemcontinuous}
\end{equation} 
where 
\begin{align*}
    \set X_C = \left\{\bm x \in \mathbb{R}^{n} \mid \bm A^\top \bm x \le \bm b, \ \bm P^\top \bm x = \bm s, \ \bm h(\bm x) \le \bm 0, \ \bm 0 \leq \bm x_B \leq \bm 1  \right\}.
\end{align*}

\medskip\noindent
{\bf Step 2: Convex relaxation via RPT.} RPT consists of two steps: a reformulation step and a perspectification step. In the reformulation step, we pairwise multiply all constraints with each other to obtain additional redundant constraints, see \cite{bertsimascone} for a guide on how to multiply any two conic constraints. In the perspectification step, we convexify the nonconvex components in both the objective as well as the additional redundant constraints by reformulating them in their perspective form and linearizing the product terms. To be more specific, let $f_k$ be an SLC function as given by \eqref{eq:SLCfunction}, that satisfies Assumption~\ref{assumption:SCC}. \cite{bertsimas2023rptbb} show that one can obtain the following convexification of $f_k$:
	\begin{equation} \sum_{i \in \set I_k} (q_{ik} - \bm d_{ik}^\top \bm x) c_{ik}\left(\frac{ q_{ik} \bm x -  \bm X \bm d_{ik} }{q_{ik} - \bm d_{ik}^\top \bm x}	\right), 
    \label{pf}
 \end{equation}
by first multiplying and dividing the argument of~$c_{ik}$ by $\left(q_{ik} - \bm d_{ik}^\top \bm x \right)$ for every $i \in \set I_k$ and then linearizing the quadratic terms $\bm x \bm x^\top$ in the argument of the reformulated $f$ with $\bm X \in \mathbb S^{n_x}$. 

RPT is then applied to \eqref{eq:genericproblemcontinuous} to obtain the following convex relaxation:
\begin{equation} 
 \label{eq:genericconvexCL} 
		\begin{array}{cll}
			\displaystyle \min_{\bm x, \bm \tau, \bm X, \bm V} & \displaystyle \quad   \tau_0 + \sum_{i \in \set I_k} (q_{i0} - \bm d_{i0}^\top \bm x) c_{i0}\left(\frac{ q_{i0} \bm x -  \bm X \bm d_{i0} }{q_{i0} - \bm d_{i0}^\top \bm x}\right) \\
			{\rm s.t.} & \displaystyle \quad  \tau_k + \sum_{i \in \set I_k} (q_{ik} - \bm d_{ik}^\top \bm x) c_{ik}\left(\frac{ q_{ik} \bm x -  \bm X \bm d_{ik} }{q_{ik} - \bm d_{ik}^\top \bm x}\right)  \le 0, & \quad k \in \set K,\\
			&\quad \displaystyle \bm b \bm x^\top \bm A + \bm A^\top \bm x \bm b^\top  \le  \bm A^\top \bm X \bm A + \bm b \bm b^\top, \\
            &\displaystyle \quad s_{\ell} \bm x - \bm X \bm P_{\ell}  = \bm 0, &\quad \ell \in \{1, \ldots, m_2\}, \\
            &\displaystyle \quad s_{\ell} \bm \tau - \bm V \bm P_{\ell}  = \bm 0, &\quad \ell \in \{1, \ldots, m_2\}, \\
			&\displaystyle \quad (b_\ell - \bm A_\ell^\top \bm x) \bm  h \left( \frac{b_\ell \bm x  - \bm X \bm A_\ell }{b_\ell - \bm A_\ell^\top \bm x} \right) \le \bm  0, & \quad \ell \in \{1, \ldots, m_1\}, \\
			& \displaystyle \quad  (b_\ell - \bm A_\ell^\top \bm x) \bm  c_0 \left( \frac{b_\ell \bm x - \bm X \bm A_\ell}{b_\ell - \bm A_\ell^\top \bm x} \right) \le b_\ell\bm \tau - \bm V \bm A_\ell,   & \quad \ell \in \{1, \ldots, m_1\}, \\
   &\quad \displaystyle X_{ii} = x_i &\quad i \in \set B, \\
   &\displaystyle \quad X_{ii} \geq 0, &\quad i \in \{B + 1, \ldots, n_x\}, \\
			& \displaystyle \quad (\bm x, \bm \tau) \in \set T, \\
               &\displaystyle \quad X_{ij} \leq x_i &\quad i > j \in \set B, \\
   &\displaystyle \quad X_{ij} \leq x_j &\quad i > j \in \set B, \\
   &\displaystyle \quad X_{ij} \geq x_i + x_j - 1 &\quad i > j \in \set B, 
		\end{array}
	\end{equation}
where $\set T = \{(\bm x, \bm \tau) \in \R^{n_x} \times \R^{K+1} \mid \bm x \in \set X_C, \ \bm c_0(\bm x) \leq \bm \tau \}$, and $\bm c_{0} (\bm x) = [c_{00}(\bm x) \ c_{01}(\bm x)  \ \cdots \  c_{0K} (\bm x)]^\top \subseteq (-\infty, + \infty ]^{K+1}$. Here $\bm \tau$ denotes the epigraphical variables for the convex components in the nonconvex SLC functions of~\eqref{eq:genericproblemcontinuous}. Multiplying these extra epigraphical constraints with the existing convex constraints can result in a tighter convex relaxation \cite[Theorem 1]{bertsimas2023rptbb}. Observe that, when RPT is applied in the presence of relaxed binary variables, the well-known McCormick inequalities (given by the last three inequalities in \eqref{eq:genericconvexCL})  for binary variables naturally emerge within the convex relaxation. Moreover, as $x_i^2$ remains the same value for binary variables, we can add the constraint $X_{ii} = x_i$ for all binary variables $x_i$.

\begin{remark}
    One might expect that the perspective function $tf(\bm x/t)$ poses numerical challenges, especially at $t=0$. However, this is for almost all practical optimization problems no issue. The reason lies in a key property of perspective functions: if the epigraph of a convex function $f(\bm x)$ is conic representable, then also the epigraph of its perspective function (in the same cone). Since nearly all practical problems can be modeled using the five basic cones (linear, conic quadratic, power, exponential, semidefinite), and efficient solvers for these cones are now widely available, perspective functions can be handled effectively.
\end{remark}

\medskip\noindent
{\bf Step 3: Branch-and-Bound framework.} We propose a hybrid branching strategy that combines the eigenvector-based branching of \cite{anstreicher2022solving} used in RPT-BB for continuous problems with classical branching techniques from mixed-integer programming (MIP). Our hybrid branching strategy works as follows. We start by solving the RPT relaxation at the root node $N_0$, and obtain a lower bound $\text{LB}_0$ and an optimal solution $(\bm x^*, \bm X^*)$. If this solution satisfies all integrality constraints and yields a matching upper bound, the algorithm terminates. If the solution does not satisfy all integrality constraints, that is, there exists at least one $i \in \mathcal{B} \cup \mathcal{D}$, such that $x^*_i \notin \mathbb{Z}$, or if the lower bound does not match the upper bound, we proceed with the branching. 

Branching involves two interdependent choices: first, whether to branch on the integer domain (using one or more integrality constraints) or on the continuous domain; and second, if we branch on the integer domain, which branching rule to apply. At each node $N$, we decide between integer branching and continuous branching by checking how ``fractional" the current RPT solution is on the integer variables. More specifically, we compute the maximum integrality violation of the current RPT solution:
\begin{align*}
    \Delta = \max_{i \in \set B \cup \set D} |x_i^* - x_i^{\text{int}}|,
\end{align*}
where $x_i^{\text{int}} = \lceil x_i^* + 0.5 \rceil$. If $\Delta$ is greater than the fractionality threshold $\delta$ we perform integer branching, otherwise we branch on the continuous variables. 

The MIP literature offers a variety of integer branching rules, including most or least fractional branching, pseudocost branching, strong branching, reliability branching, randomized selection, hybrid branching, and more (see \cite{Benichou71}, \cite{achterberg05}, \cite{achterberg07}, \cite{achterberg09}). Our approach uses most fractional branching, which selects the integer variable whose value in the current LP relaxation is farthest from an integer, aiming to quickly reduce fractional infeasibility in the branch-and-bound tree.
However, users may substitute any integer-branching rule that best exploits the structure of their particular nonconvex mixed-integer problem. 

If $\Delta \leq \delta$, we branch in the continuous domain using eigenvector branching. 
That is, we select $\bm X^*_{C}$ to be the sub-matrix of $\bm X^*$ that only contains the rows and columns of $\bm X^*$ that correspond to the continuous variables. 
Now, we find the eigenvector $\bm v$ of $\bm X^*_{C} - \bm x^*_{C}\bm x^{*\top}_{C}$ with the largest eigenvalue $\lambda$. The branching hyperplane is then $\set H = \{\bm x \mid \bm f^\top \bm x = \lambda \}$, where 
$\bm f = \begin{bmatrix} (\bm 0^{|\set B \cup \set D|})^\top, \bm v^\top \end{bmatrix}^\top$. Hence, we form two child nodes $N_{\ell}$ and $N_{r}$ and in each node we solve the original problem with the new constraint sets:
\begin{align*}
    & \mathcal{X}_{N_{\ell}} = \{\bm x\in \mathcal{X}_{N} | \bm f^\top \bm x \leq \lambda\}, \\
    & \mathcal{X}_{N_r} = \{\bm x\in \mathcal{X}_{N} | \bm f^\top \bm x \geq \lambda\},
\end{align*}
where $\mathcal{X}_{N}$ denotes the feasible region of the RPT relaxation at node $N$. The proposed branching strategy is summarized in Algorithm \ref{alg:branching_strategy}.

\small{
 \begin{algorithm}[h]
        \caption{Branching strategy \label{alg:branching_strategy}} 
        \textbf{Input}: ($N, \set X_N$). \\ 
        \textbf{Output}: $(N_{\ell},N_r)$.
      \begin{algorithmic}[1]
            \STATE Solve RPT relaxation at node $N$ and obtain optimal solution $(\bm x^*, \bm X^*)$ 
            \IF{$\exists i \in \set B \cup \set D$ such that $x_i^* \notin \mathbb{Z}$}
            \STATE Set $\Delta = \max_{i \in \set B \cup \set D} |x_i^* - x_i^{\text{int}}|$
            \ELSE
            \STATE Set $\Delta = 0$
            \ENDIF
            \IF{$\Delta > \delta$}
            \STATE Perform integer branching
            \ELSE
            \STATE Take $\bm v$ as the eigenvector corresponding to the largest eigenvalue of $\bm X_{C}^* - \bm x_{C}^* (\bm x_{C}^*)^\top$ 
            \STATE Take $\lambda = \bm v^\top \bm x_{C}^*$ 
            \STATE Set $\bm f = \begin{bmatrix} (\bm{0}^{|\mathcal{B} \cup \set D|})^\top, \bm v^{\top} \end{bmatrix}^\top$
            \STATE Create child node $N_{\ell}$ with feasible region $\set X_{N_{\ell}} = \left\{ \bm x \in \set X_N \mid \bm f^\top \bm x \le \lambda \right\}$ 
            \STATE Create child node $N_r$ with feasible region $\set X_{N_r} = \left\{ \bm x \in \set X_N \mid \bm f^\top \bm x \ge \lambda \right\}$
            \ENDIF
            \RETURN $N_{\ell},N_{r}$
         \end{algorithmic}
    \end{algorithm}}
\normalsize 

After creating the two child nodes, we solve the RPT relaxation on each node, in which the new inequality is multiplied with each already existing constraint, to obtain their respective lower and upper bounds. We then compare the lower bounds: if the smaller of the two lower bounds matches the best upper bound obtained so far, the algorithm terminates with an optimal solution. Otherwise, we add the child nodes to a priority queue and select the next node to branch on as the one with the smallest lower bound. 

\begin{remark}
    When branching on an integer variable, for example by imposing $x_1 \leq 5$ or $x_1 \geq 6$, in a similar way we can create two child nodes $\set X_{N_1} = \{ \bm x \in \set X_N \mid x_1 \leq 5\}$ and $\set X_{N_2} = \{ \bm x \in \set X_N \mid x_1 \geq 6\}$, and solve the RPT relaxation on each node in which the new inequality is multiplied with each already existing constraint.
\end{remark}

If we are dealing with a problem that has only discrete variables and we find a solution $\bm x^*$ such that all variables are integer but the lower and upper bounds do not match, we can no longer do eigenvector branching. In this case, we skip the current node and look at the next node in the priority queue. The RPT-BB algorithm for nonconvex mixed integer optimization problems is summarized in Algorithm \ref{alg:rpt_bb}.

    \begin{algorithm}[h]
        \caption{Branch and bound via RPT \label{alg:rpt_bb}} 
        \textbf{Input}: ($N_0, \text{Lb}^0$, $\text{Ub}^0$, $\varepsilon$). \\ 
        \textbf{Output}: $(\bm x^*, \text{Lb}, \text{Ub})$.
      \begin{algorithmic}[1]
            \STATE $\text{Lb} \gets \text{Lb}^0 $
            \STATE $\text{Ub} \gets \text{Ub}^0 $
            \STATE ACTIVE $\gets \{N_0\}$
            \WHILE{$\text{Ub} - \text{Lb} > \varepsilon$}
            \STATE $j \gets \argmin_{i \in \text{ACTIVE}} \text{Lb}^i$ 
            \STATE Partition node $N_j$ into two child nodes $N_{j_1},N_{j_2}$ by applying Algorithm \ref{alg:branching_strategy}
            \FOR{$i = 1,2$}
            \STATE Solve the RPT relaxation corresponding to $N_{j_i}$ and obtain Lb$^{j_i}$ and Ub$^{j_i}$. 
            \ENDFOR
            \STATE 
            $\text{Ub} \gets \min\{\text{Ub}^j, \text{Ub}^{j_1}, \text{Ub}^{j_2}\}$ 
        \FOR{$i = 1,2$}
                \IF{$\text{Lb}^{j_i} < \text{Ub}$}
                \STATE ACTIVE $\gets$ ACTIVE $\cup \{j_i\}$
                \ENDIF
                \ENDFOR
        \STATE Lb $\gets \min\{\text{Lb}^{j_1}, \text{Lb}^{j_2}\}$
        \STATE 
         $\text{ACTIVE} \gets \text{ACTIVE} \setminus \{j\}$
         \ENDWHILE
         \end{algorithmic}
    \end{algorithm}



\section{Numerical Experiments} \label{sec:numexp}
We evaluate the performance of our approach on five nonconvex mixed-integer optimization problems, taken from \cite{bertsimas2023rptbb} . Numerical experiments are performed on an Apple M3 Pro chip. All computations for RPT-BB and SCIP are conducted with MOSEK version 10.2.0, Gurobi version 11.0.2 \citep{Gurobi}, SCIP version 0.11.6 \citep{scip}, and implemented using Julia 1.11.3 and the Julia package \texttt{JuMP.jl} version 1.23.2. All computations for BARON are conducted with BARON version 2021.1.13 \citep{baron} implemented using the Python package \texttt{pyomo} version 6.8.0. Inside RPT-BB, we use Gurobi for the linear optimization problems and MOSEK for the nonlinear optimization problems. In all branch and bound implementations the optimality gap is set to $10^{-4}$. Moreover, we take the branching threshold $\delta = 0$ (see Section \ref{sec:branch_strat} for a comparison of branching strategies), as this choice consistently led to the best performance across both the mixed-binary and mixed-integer domains in our experiments. We compare our approach, when applied with (RPT-SDP-BB) and without the LMI (RPT-BB), with BARON, and with SCIP.

For each numerical experiment, we report the best objective value found, the computation time in seconds, the number of decision variables, the number of integer branchings, and the number of eigenvector branchings, denoted by \textbf{Opt}, \textbf{Time}, \textbf{Dims}, \textbf{Int Hyp}, and \textbf{Eigen Hyp}, respectively. Additionally, all solvers are given a time limit of one hour. If a solver exceeds this limit without proving optimality, we record the running time as $\bm{3600^*}$; if it fails to produce any feasible solution, we set the optimal value to $\bm{\infty}$ and if a problem is provably infeasible, we denote the optimal value as \textbf{Infeasible}. The results in Tables \ref{tab:numexp_somol}, \ref{tab:numexp_lse}, and \ref{tab:numexp_lmp}, are averaged over 5 randomly generated instances. When solving problems in the mixed-binary and mixed-integer domains, we consider the case where the first $\lfloor{\frac{n}{2}}\rfloor$ variables are integer, and the remaining $n -\lfloor{\frac{n}{2}}\rfloor$ are continuous.

\subsection{Sum-of-max-linear-terms maximization}
We consider the sum-of-max-linear-terms (SOMOL) maximization problem from \cite{bertsimas2023rptbb} 
\begin{align} \label{eq:smlt} 
		\begin{array}{cll}
			\displaystyle \max_{\bm x \in \mathcal{X}} &\quad  \displaystyle  \sum_{\ell \in \mathcal{L}} \max_{j \in \mathcal{J}_{\ell}} \{\bm{A}_{j}^\top\bm{x} + b_j\},
		\end{array}
\end{align}	
\noindent where $\bm A \in \R^{n_x \times n_y}, \bm b\in \R^{n_y}$. We consider two sets of constraints, defined as follows
\begin{align*}
    & \mathcal{X}_1 = \left\{\bm x \equiv (\bm x_B, \bm x_D, \bm x_C) \ \middle| \ \bm D^\top \bm x \leq \bm d\right\}\text{,} \\ 
    & \mathcal{X}_2 = \left\{ \bm x\in\mathcal{X}_1 \ \middle| \ \log\left(\sum_{i=1}^{n_x}\operatorname{exp}(x_i)\right)\leq \alpha\right\}\text{,}
\end{align*}
where $\bm D \in \R^{n_x \times m}, \bm d \in \R^{m}$. 

We examine five different problem instances of varying complexity, see  \citet[Appendix D.1]{bertsimas2023rptbb}. We run each problem instance for each of the four variable domains, i.e., binary, mixed binary, integer, and mixed integer, and record the results in Table \ref{tab:numexp_somol}.

\begin{table}[h]
\resizebox{\columnwidth}{!}{

 \begin{tabular}{cccccccccccccccccc}
\toprule 
\multirow{2}{*}{} & \multirow{2}{*}{$\set X$} & \multirow{2}{*}{\#} & \multirow{2}{*}{$n_x$} &\multicolumn{4}{c}{\textbf{RPT-BB}}&\multicolumn{4}{c}{\textbf{RPT-SDP-BB}}
&\multicolumn{2}{c}{\textbf{BARON}}&\multicolumn{2}{c}{\textbf{SCIP}}
\\ \cmidrule(r){5-8}\cmidrule(r){9-12}\cmidrule(r){13-14}\cmidrule(r){15-16}
&&&&Opt& Time & Int Hyp  & Eigen Hyp&Opt& Time & Int Hyp  & Eigen Hyp & Opt & Time & Opt & Time \\
\bottomrule

\multirow{10}{*}{B} & \multirow{5}{*}{$\set X_1$}  & 1 & 5& 9.85 &4.66&	0 & 0 & 9.85 & 5.72	& 0	& 0	& 9.85 &	\textbf{0.05}	&	9.85	&2.52 \\

&& 2 & 5 & 83.33	& \textbf{0.02}	&0	& 0	&83.33&	 0.37&	0	& 0 &83.33	&0.18		&83.33	&0.03\\

&& 3 & 20 & 208.80	& 52.41	&182 & 0 &	208.80	&32.46	&0 & 0 &	208.80	&362.38	&	208.8&	\textbf{13.67}  \\

&& 4 & 10 &32.82	& \textbf{0.01}	& 0	& 0 &	32.82 &	0.03 &	0 & 0	&	32.82&	0.18	&	32.82 & 2.41 \\

&& 5 & 20 & 791.37 &	\textbf{1.30} &	0 & 0	&	791.37 &	32.47&	0	& 0 &791.37 &	10.87	&	791.37	& 1.84\\

\cline{2-16}

 &\multirow{5}{*}{$\set X_2$}  & 1 & 5 & 9.85	& 5.81 &	0 & 0	&	9.85 &	6.81 &	0 & 0	&	9.85 &	\textbf{0.05}	&	9.85 &	2.77 \\
 
&& 2  & 5 & 83.33	& \textbf{0.03} &	0 & 0	&	83.33	& 0.57	& 0 & 0	&	83.33 & 	0.17	&	83.33 &	2.85\\

&& 3 & 20 & 208.80	& 155.27	& 162 & 0	 &	208.80 &	401.81 &	0 & 0	&	208.80	& 464.99	&	208.80 &	\textbf{11.05} \\

&& 4  & 10 & 27.96	& 1.51 &	20 & 0	&	27.96 &	9.03 &	19	 & 0 &	27.96 &	\textbf{0.56}	&	27.95 &	2.71   \\

&& 5 & 20  & 791.37 &	\textbf{2.04} & 	0 & 0	&	791.37	& 584.54 &	0	& 0 &	791.37	& 12.62	&	791.37 &	3.38 \\
\bottomrule


\multirow{10}{*}{MB}& \multirow{5}{*}{$\set X_1$} & 1 &  $5$  & 10.93 &	8.80	&0&	0&		10.93&	10.21&	0	&0	&	10.93	& \textbf{0.05}	&	10.93&	2.47\\

&& 2 & $5$ & 89.52	&\textbf{0.02}	&0&	0	&	89.52	&0.38	&0	&0		&89.52 &	0.19	&	89.52	&0.05\\

&& 3 & $20$ & 244.51	& 73.30	&381 &	12	&	244.51	&33.19&	0	&0		&244.51	&3600*	&	244.51	&\textbf{20.91}\\

&& 4 &$10$ & 68.37 &\textbf{0.01} &	0 &	0	&	68.37	& 0.03 &	0	&0		&68.37	&1.38	&	68.37&	0.05\\

&& 5  &$20$ & 2649.33 &	\textbf{40.21}  &	92	& 31 &	2649.33	&357.17 &	5	&1		&2649.33	&3600*		&2561.98	&$3600^*$\\

\cline{2-16}

&\multirow{5}{*}{$\set X_2$}& 1 & $5$ & 10.93	&8.30 &	0 &	0	&	10.93	&10.38 &	0	&0		&10.93 &	\textbf{0.05}	&	10.93	&2.87  \\

&& 2 & $5$ & 89.52	&\textbf{0.03}&	0	&0		&89.52	&0.74&	0	&0		&89.52 &	0.22	&	89.52&	0.26  \\

&& 3 &	20 &244.51	&288.71 &	377&	12	&	244.51&	88.66 &0	&0	&	244.51&	3600* &	244.51 &	\textbf{28.41}	 \\

&& 4 &	10 & 31.18	&0.63&	7	&0		&31.18	&1.70&	3	&0	&	31.18	&1.69	&	31.18	&\textbf{0.18}	 \\

&& 5 &	20 & 1483.24&	\textbf{155.84}&	71	&47	&	1483.24	&3600*	&18&	2&		1483.24	&3600*	&	1483.24	&3600*	\\
\bottomrule


\multirow{10}{*}{I} &\multirow{5}{*}{$\set X_1$}  & 1 & 	5 & 21.78&	4.66	&2	& 0 &	21.78	&5.63&	2& 0	&	21.78	&\textbf{0.05}	&	21.78&	2.67	\\

&& 2 &	5 &219.28&	0.07	&2 & 0	&	219.28	&1.66&	2 & 0 &	219.28	&0.11	&	219.28	&\textbf{0.03}	\\

&& 3 & 20 &1051.38 &	 547.62 &	2731 & 0 & 1051.38	&711.43	&19	 & 0 	&1051.38 &3600* &	1051.38	&\textbf{6.12}  \\

&& 4 &	10 &113.32	& 0.09&	5 & 0	&	113.32	&0.44&	5 & 0	&	113.32 &	2.83	&	113.32	&\textbf{0.04}	\\

&& 5  & 20 & 2954.89 &	\textbf{21.95} &	67 & 0	&2954.89	& 1855.28	&49& 0	&	2954.89	&3600*	 &	2603.73	&3600*	\\

\cline{2-16}

& \multirow{5}{*}{$\set X_2$}  & 1  &	5 & 12.61	& 6.35	&5	& 0 &	12.61	& 7.40 &	5	& 0 &	12.61&	\textbf{0.11}	&	12.61	& 2.87	\\

& & 2 &	5 & 129.25	& 0.45	&5	& 0 &	129.25	&5.83	&4	& 0 &	129.25	&0.20	&	129.25&	\textbf{0.03}	\\

&& 3 &	20 & 793.87	&2709.43	&2774 & 0 &	793.87	&663.88&	5	& 0	&793.87	& 3600*	&793.87	&\textbf{12.65}	\\

&& 4 &	10 & 31.24	&2.23&	28	&  0 &	31.24	&10.06	&18 & 0	&31.24&	2.00	&	31.24 &	\textbf{0.13}	 \\

&& 5 &	20 &1562.71&	\textbf{3600*} &2046 & 0	&1562.59	&3600*	&22  & 0 &	1552.46 &	3600*	&	1549.95	&3600*	 \\
\bottomrule


\multirow{10}{*}{MI} & \multirow{5}{*}{$\set X_1$}& 1 &  $5$ & 22.85	&8.19&	1&	0	&	22.85	&9.31&	1&	0	&	22.85	& \textbf{0.05}	&	22.85&	2.55 \\

&& 2 & $5$ & 224.93	&\textbf{0.04}&	1&	0	&	224.93	&0.96&	1	&0	&	224.93	&0.14	&	224.93	&\textbf{0.04}\\

&& 3 & $20$ & 1069.34	& 47.78	&218	&0	&	1069.34	& 87.01&2&	0	&	1069.34	&681.44	&	1069.34	&\textbf{1.59}\\

&& 4 &$10$ & 113.32&	0.09	&5&	0		&113.32	&0.26&	5&	0	&	113.32	&1.04	&	113.32	&\textbf{0.06}\\

&& 5 &$20$ & 2987.78	& \textbf{5.17} &12	&0	&	2987.78	&215.98&	5	&0		&2851.38	&3600*		&2603.13	&3600*\\

\cline{2-16}
&\multirow{5}{*}{$\set X_2$}  & 1& $5$  & 13.69	&9.07&	3&	0	&	13.69&	11.48&	3&	0	&	13.69	&\textbf{0.18}	&	13.69&	2.71 \\

&& 2& $5$ & 134.48	&0.23 &	2&	0	&	134.48	& 3.88 &	2	&0	&	134.48	&0.24	&	134.48	&\textbf{0.06}  \\

&& 3 &	20 &809.78&	94.31	&125	&0	&	809.78	&438.38 &3&	0	&	809.78	&848.49	&	809.78	&\textbf{16.56}	 \\

&& 4 &	10 & 31.91	&1.05&	11	&0	&	31.91	& 5.20 &	9&	0	&	31.91	&1.96	&	31.91	&\textbf{0.21}	 \\

&& 5 &	20 & 1593.48&	\textbf{3600*}	&2406	&0	&	1593.48&	\textbf{3600*} &	15&	0	&	1593.04	&3600*	&	1569.19	& 3600*	\\
\bottomrule

\end{tabular}
}
\caption{Comparison of the approaches RPT-BB, RPT-SDP-BB, BARON, and SCIP for the binary (B), mixed-binary (MB), integer (I), and mixed-integer (MI) domains for Problem \eqref{eq:smlt}. \label{tab:numexp_somol}}
\end{table} 

From Table \ref{tab:numexp_somol}, we observe that while all methods are able to solve problem instances 1, 2, and 4, in less than a minute across all domains, BARON outperforms all methods for problem instance 1, while for problem instances 2 and 4 all methods have similar performance.
For problem instance 3, all approaches are able to find the global optimum within the time limit, except for BARON in the integer and mixed integer case, but SCIP is the most efficient method. RPT-BB is the most efficient method for problem instance 5, for both $\set X_1$ and $\set X_2$ across all variable domains. Moreover, for problem instance 5, we observe that BARON and SCIP are not able to prove optimality within the time limit in the mixed binary, integer, and mixed integer case. For the integer and mixed integer cases, we observe that all methods are not able to find the optimum within the time limit for $\set X = \set X_2$.



\subsection{Quadratic constraint quadratic optimization} \label{sec:qcqp}
We next consider the quadratic constraint quadratic optimization problem (QCQP) from \cite{bertsimas2023rptbb} 
\begin{equation} 
\begin{aligned}
\label{eq:QCQPnconvex}
& \displaystyle \underset{\bm x \in \mathcal{X}_1}{\text{max}}
& & \displaystyle \bm x^\top \bm P_0 \bm x +  \bm q_0^\top \bm x + r_0 \\
& {\rm s.t.}
& & \displaystyle \bm x^\top \bm P_k \bm x + \bm q_k^\top \bm x + r_k \leq 0, &\quad k \in \set K,
\end{aligned}  
\end{equation}  
where $\bm P_k \in \R^{n_x \times n_x}, \bm q_k \in \R^{n_x}$. The matrices $P_k$ are not necessarily positive semidefinite and, hence, these constraints are not necessarily convex.

We examine ten different problem instances of varying complexity, see  \cite[Appendix D.4]{bertsimas2023rptbb}. Problems 1-5 correspond to minimizing a quadratic function subject to convex inequalities only, that is, $\mathcal{K} =\emptyset$. Problems 5-10 correspond to minimizing a quadratic function subject to both convex and nonconvex inequalities. We run each problem instance in each of the four variable domains and record the results in Table \ref{tab:numexp_qcqp}.

\begin{table}[h]
\resizebox{\columnwidth}{!}{
 \begin{tabular}{ccccccccccccccccc}
\toprule 
 \multirow{2}{*}{$\set X$} & \multirow{2}{*}{\#}&  \multirow{2}{*}{$n_x$} & \multicolumn{4}{c}{\textbf{RPT-BB}}&\multicolumn{4}{c}{\textbf{RPT-SDP-BB}}
&\multicolumn{2}{c}{\textbf{BARON}}&\multicolumn{2}{c}{\textbf{SCIP}}
\\ \cmidrule(r){4-7}\cmidrule(r){8-11}\cmidrule(r){12-13}\cmidrule(r){14-15}
&&& Opt& Time & Int Hyp  & Eigen Hyp&Opt& Time & Int Hyp  & Eigen Hyp & Opt & Time & Opt & Time   \\
\bottomrule

\multirow{10}{*}{B}& 1 & 20 &35.50 &      8.50   &  16        & 0&	35.50	& 10.68	&16	&0&	35.50	& 0.12		&35.50 & \textbf{2.56}\\

& 2& 20 & 360.00  &	0.39 &	0&0	&	360.00	& 0.39 &	0& 0	&	360.00	& 0.05		& 360.00	& \textbf{0.05} \\

& 3 & 10 & 275.55 &	\textbf{0.05}  &	0 & 0	&	275.55 &	0.07	& 0	& 0&	275.55	& \textbf{0.05}	&	275.55	& \textbf{0.05} \\

& 4 & 50 &32033.02	& 0.51 &	0 &0	&	32033.02	& 0.78 &	0	& 0 &	32033.02&	0.17	&	32033.02 &	\textbf{0.01} \\

& 5 & 100 &  253526.05	& 5.36 &	0 &0	&	253526.05	& 14.11 &	0	&0	& 253526.05	& 0.79	&	253526.05	& \textbf{0.03}\\

 & 6 & 8 &	Infeasible & 0.33  &	0	&0 &	Infeasible &	0.37 &	0  &0 &	 Infeasible&	 0.05	&	Infeasible &	0.09\\
 
& 7& 12 & Infeasible	&0.01 &	0 & 0 &	Infeasible	& 	0.01 & 0	&0&	Infeasible & 	0.05	&	Infeasible &	0.01 \\

& 8 & 16 & Infeasible & 0.01	&  0 & 0 &	Infeasible & 0.02  &	0 & 0 	&	Infeasible & 0.13	&	Infeasible &	0.01\\

& 9 & 30 & 58.00 &	38.61 & 383	 &0 &  58.00  & 79.11	& 296 &0 &	58.00	& 6.82		& 58.00	& \textbf{2.66}  \\

& 10 & 40 & Infeasible & 1.45 & 6 & 0 &	Infeasible	& 0.48 &	0 &0 &	Infeasible	& 1.85 &	Infeasible &	 0.42 \\
\hline

\multirow{10}{*}{MB} & 1 & 20 & 247.76  &	7.68  &	2 &	3	&	247.76	&9.16&	1	& 2&		247.76 &	0.26	&	247.76	& \textbf{0.25}  \\

& 2& 20 & 737.76 &	0.85 &	2	& 3	&	737.76	& 1.10 &	1	& 2 &		737.76 &	\textbf{0.24}	&	737.76 &	0.29	\\

& 3 & 10 & 2706.69	& 0.06	& 0	& 0	&	2706.69 &	0.06 &	0 &	0	&	2706.69	& \textbf{0.03}	&	2706.69	&0.05	\\

& 4 & 50 &  77333.33	&\textbf{3.93} &	3&	0	&	77333.33	& 18.79 &	4	&0		&77333.33	& 46.94 &	77333.33	& 19.70	\\

& 5 & 100 & 572687.27	&6.40 &	0&	0	&	572687.27	& 17.87	&0&	0		&572687.27	&0.77		&572687.27	&\textbf{0.04}	\\

& 6 & 8 & 408.00	& 0.30	& 0	& 0	&	408.00&	0.31 &	0	& 0	&	408.00	& \textbf{0.03}	&	408.00	& 0.10	\\

& 7& 12 &  282.06&	0.65	&9&	0		&282.06&	0.55&	5&	0		&282.06 &	0.12	&282.06&	\textbf{0.07}	 \\

& 8 & 16 & 521.35&	0.21 &	6&	0		&521.35	&0.29&	3&	0	&	521.35&	\textbf{0.18}	&	521.35&	0.25	\\

& 9 & 30 & 2149.32&	4.98	&38&	0	&	2149.32	& 5.03	&14&	0&		2149.32	&1.35		&2149.32	&\textbf{0.67}	\\

& 10 & 40 & 2421.69	&174.83	&506&	0		&2421.69	&31.04&	35	&0	&	2421.69&	17.42		&2421.69	&\textbf{5.26}\\
\hline

\multirow{10}{*}{I} & 1 & 20 & 370	&8.64&	30 & 0	&	370	&12.96&	30	& 0 &	370	&1.91		&370	&\textbf{2.69}	\\

& 2& 20 & 860	&2.22&	32 & 0	&	860	& 5.18	&32	 & 0 &	860	&1.63		&860	& \textbf{0.89} \\

& 3 & 10 &  6888.78	&0.05	&0	& 0 &	6888.78	&0.05	&0 & 0	&	6888.78	&\textbf{0.03}	&	6888.78	&0.05 \\

& 4 & 50 & 97927.97	& 70.10	&100 & 0	&	97927.97& 221.15 &100	& 0 &	95625.73&	$3600^*$	&	97927.97&	\textbf{43.02} \\

& 5 & 100 & 772790.63	&2794.04	&132  & 0	&	772790.63&	$3600^*$	&76	& 0 &	764805.87& $3600^*$	&	772790.63&	\textbf{298.9}\\

 & 6 & 8 & 	1716	&0.37	&1	& 0 &	1716	& 0.38	&1	& 0 &	1716	&\textbf{0.08}	&	1716	&0.12 \\
 
& 7& 12 & 4482	&0.14&	14 & 0	&	4482&	0.43 &	15 & 0	&	4482&	0.14	&	4482	&\textbf{0.06} \\

& 8 & 16 & 17045	&\textbf{0.11}	&8	& 0 &	17045&	0.41 &	8 & 0	&	17045&	0.19		&17045	& 0.12 \\

& 9 & 30 & 24170.0	&$3600^*$	&23003   &   11008&	46558 &	\textbf{152.63}	& 588 & 0	&	46558&	$3600^*$	&	41217	& $3600^*$   \\

& 10 & 40 & Inf&	$3600^*$	&2208	& 3104 &	44574.0	&$3600^*$	& 5448	& 0 &	51542&	\bm{$3600^*$}	& -1589	& $3600^*$ \\
\hline

\multirow{10}{*}{MI} & 1 & 20 & 375.84&	7.87 & 4&	0	&	375.84	& 8.90&	4	&0	&	375.84	& 0.42	&	375.84	&\textbf{2.94}	 \\

& 2& 20 & 865.84	& 0.79	&4	&0	&	865.84	& 1.16	&4&	0	&	865.84&	0.36	&	865.84	&\textbf{0.21} \\

& 3 & 10 & 6888.78	&0.06&	0	&0	&	6888.78&	0.07	&0	&0&		6888.78	& \textbf{0.03}	&	6888.78 &	0.05	\\

& 4 & 50 & 98331.2& \textbf{2.43} &	2 &	0	&	98331.2	&6.95 &2	&0		&98050.57&	$3600^*$	&	98331.2	& 55.73	\\

& 5 & 100 & 774403.58	&  \textbf{90.64} &	2 &	0	&	774403.58	& 138.70	& 2	& 0		& 774403.58 &	$3600^*$	&	774403.58 &	590.53\\

& 6 & 8 & 1717.80 &	0.41 &	0 &	0	&	1717.80	&0.42 &	0	&0		&1717.8	&\textbf{0.05}	&	1717.80 &	0.14\\

& 7& 12 & 4505.24	& 0.07 &	2 &	0	&	4505.24	&0.11 &	2	&0		&4505.24	&0.11	&	4505.24	&\textbf{0.04}	 \\

& 8 & 16 & 17063.66 &	0.16 &	3 &	0	&	17063.66	& 0.24 &	3&	0		&17063.66&	0.19	&	17063.66	&\textbf{0.09}\\

& 9 & 30 & 42821.70	&$3600^*$	&9029	&3151	&	46619.96&	\textbf{202.84}	& 551 &	0	&	46619.96 &	3515.34	 &	46619.99	& 579.74	\\

& 10 & 40 & 23914.80	&$3600^*$&	1482&	1185	&	52210.15&	\bm{$3600^*$} &	3729	&0	&	51816.65	& $3600^*$	&	52210.15	& \bm{$3600^*$}	\\
\hline
\end{tabular}
}
\caption{Results for Problem \eqref{eq:QCQPnconvex} when applying RPT-BB, RPT-SDP-BB, BARON, and SCIP the binary (B), mixed-binary (MB), integer (I), and mixed-integer (MI) domains. \label{tab:numexp_qcqp}}
\end{table}

From Table \ref{tab:numexp_qcqp}, we observe that in all variable domains, for the simpler problem instances 1, 2, 3, 6, 7, and 8, with fewer decision variables, all solvers are very efficient, solving the problem in a few seconds whenever feasible. Additionally, when we optimize over only binary variables, we see that problems 6, 7, and 8 are infeasible, with both RPT-BB and RPT-SDP-BB detecting infeasibility at the root node. 

For problem instance 4, SCIP performs significantly better than all other approaches in the integer case, for which both BARON and RPT-SDP-BB are not able to prove global optimality within the time limit, while RPT-BB performs best for the mixed binary and mixed integer case, with BARON again being unable to find the optimal solution within the time limit in the mixed integer case.

For problem instance 5, SCIP is the most efficient method in the binary, mixed binary, and integer case. However, in the mixed integer case we observe that RPT-BB and RPT-SDP-BB are significantly more efficient than both BARON and SCIP. More precisely, RPT-BB terminates in less than 4 minutes, while SCIP requires more than 50 minutes, and BARON is unable to reach provable optimality within the time limit.

Finally, for the more complex nonconvex constrained problems 9 and 10, we observe the following: For problem instance 9, in the binary and mixed binary cases SCIP is the most efficient, followed by BARON, and then by RPT-BB. In the integer case we observe that only RPT-SDP-BB solves the problem to provable optimality, while BARON, SCIP, and RPT-BB are unable to reach optimality within the time limit. In the mixed integer case, we see that SCIP again is the most efficient, followed by RPT-SDP-BB, and then by BARON and RPT-BB, which are unable to reach optimality within the time limit. For problem instance 10 we observe that it is infeasible in the binary case. Very interestingly, while RPT-BB requires 8 branches to detect infeasibilty, with the SDP constraint the algorithm detects it at the root node. In the mixed binary case SCIP appears to be the most efficient, followed by BARON, RPT-SDP-BB and then RPT-BB. In the integer and mixed integer cases we observe that no algorithm reaches optimality within the time limit.

Overall, SCIP is generally the most efficient approach. However, RPT-BB and RPT-SDP-BB demonstrate clear advantages in handling the more complex cases in the mixed integer domain. 



\subsection{Log-sum-exp maximization over linear constraints}
We also consider the log-sum-exp maximization problem subject to linear constraints from \cite{bertsimas2023rptbb}

\begin{align} \label{eq:lse} 
		\begin{array}{cll}
			\displaystyle \max_{\bm x\in \mathcal{X}} &\quad \displaystyle \log\left(\sum_{i = 1}^{n_x}\exp(x_i)\right),
		\end{array}
\end{align}	

\noindent where $\mathcal{X} = \{\bm x \ | \ \bm D^\top \bm x \bm \leq \bm d\}$, $\bm D \in \mathbb{R}^{n_x\times L'}$ and $\bm d \in \mathbb{R}^{L'}$. 

We examine five different problem instances from \cite{bertsimas2023rptbb}, where the dimension ranges from 10 to 50. We run each problem for the mixed binary and mixed integer domains and record the results in Table \ref{tab:numexp_lse}. 

\begin{table}[h]
\resizebox{\columnwidth}{!}{
 \begin{tabular}{ccccccccccccccccc}
\toprule 
 \multirow{2}{*}{$\set X$} & \multirow{2}{*}{\#}&  \multirow{2}{*}{$n_x$} & \multicolumn{4}{c}{\textbf{RPT-BB}}&\multicolumn{4}{c}{\textbf{RPT-SDP-BB}}
&\multicolumn{2}{c}{\textbf{BARON}}&\multicolumn{2}{c}{\textbf{SCIP}}
\\ \cmidrule(r){4-7}\cmidrule(r){8-11}\cmidrule(r){12-13}\cmidrule(r){14-15}
&&& Opt& Time & Int Hyp  & Eigen Hyp&Opt& Time & Int Hyp  & Eigen Hyp & Opt & Time & Opt & Time   \\
\bottomrule

\multirow{5}{*}{MB} & 1 & 10 & 3.48 & 1.91  & 0   &  0 & 3.48 & 2.43 & 0  & 0 & 3.48 & 0.09 & 3.48 & 0.43\\

& 2& 40 &  4.91 & 0.53 & 0 & 0 & 4.91 & 84.54 & 0 & 0 & 4.91 & 0.04 & 4.91 & 0.01\\

& 3 & 10 & 6.25 & 0.66  & 0 &   0.2	& 6.23 & 1.20 &  0 &       0.2 & 6.25 & 0.05 &  6.25 & 0.01\\

& 4 & 20 & 22.53 & 0.24 &  0   &  0 & 22.53 & 1.63 & 0 & 0 & 22.53 & 0.06 & 22.53 &  0.02\\

& 5 & 50 &  34.65 & 5.63 &  0 &  0 & 34.65 & 118.07 & 0 & 0 & 34.65 & 0.08 & 34.65  & 0.06 \\

\hline

\multirow{5}{*}{MI} & 1 & 10 & 10.01 & 2.29 &  2 & 0	& 10.01& 2.70 & 2 & 0 & 10.01 & 0.04 & 10.01  & 0.45\\

& 2& 40 & 40.0  & 1.02 & 0 & 0 & 40.0  &  28.39 & 0 & 0 & 40.00 & 0.08 & 40.0  &  0.01\\

& 3 & 10 &6.25 &  1.25 &  1 &   0.2	& 6.25 &  1.39 & 0.4  & 0 & 6.25 & 0.06 &  6.25 & 0.01\\

& 4& 20  & 23.35 & 0.75  & 0.8  & 0 & 23.35 & 2.14 & 0.4 &  0 & 22.53 & 0.09 & 23.35  & 0.02	\\

& 5 & 50 & 35.24  & 9.76 &   0.6  &  0 & 35.24 & 180.87 & 0.4  & 0 & 22.61 & 0.09 & 35.24  &  0.07\\

\bottomrule
\end{tabular}
}
\caption{Results for Problem \eqref{eq:lse} when applying RPT-BB, RPT-SDP-BB, BARON, and SCIP for mixed-binary (MB) and mixed-integer (MI) domains. \label{tab:numexp_lse}}
\end{table} 

From Table \ref{tab:numexp_lse}, we observe that SCIP is the most efficient method. RPT-BB, BARON and SCIP are able to solve all problems in less than a minute. Moreover, we notice that RPT-BB is able to find the optimum in the root note, in the mixed binary case, while it requires at most 2 hyperplanes in the mixed integer case.

\subsection{Linear multiplicative optimization}
We next consider the following linear multiplicative problem from \cite{bertsimas2023rptbb}

\begin{align} \label{eq:lmpopt} 
		\begin{array}{cll}
			\displaystyle \min_{\bm x\in\mathcal{X}} &\quad \displaystyle \prod_{i=1}^{n_y} \bm A_i^T \bm x + \bm b_i,
		\end{array}
\end{align}	

\noindent where $\bm A \in \mathbb{R}^{n_x\times n_y}$, $\bm b \in \mathbb{R}^{n_y}$, $\mathcal{X} = \{\bm x \ | \ \bm D^T \bm x \bm \leq \bm d, \bm A_i ^T\bm x +\bm b \geq 0\}$,  $\bm D\in \mathbb{R}^{n_x\times L}$, and $\bm d \in \mathbb{R}^{L}$. 

We examine five different problem instances from \cite{bertsimas2023rptbb}, where the dimension ranges from 5 to 40. We run each problem for the mixed binary and mixed integer domains and record the results in Table \ref{tab:numexp_lmp}. 

\begin{table}[h]
\resizebox{\columnwidth}{!}{
 \begin{tabular}{ccccccccccccccccc}
\toprule 
 \multirow{2}{*}{$\set X$} & \multirow{2}{*}{\#}&  \multirow{2}{*}{$n_x$} & \multicolumn{4}{c}{\textbf{RPT-BB}}&\multicolumn{4}{c}{\textbf{RPT-SDP-BB}}
&\multicolumn{2}{c}{\textbf{BARON}}&\multicolumn{2}{c}{\textbf{SCIP}}
\\ \cmidrule(r){4-7}\cmidrule(r){8-11}\cmidrule(r){12-13}\cmidrule(r){14-15}
&&& Opt& Time & Int Hyp  & Eigen Hyp&Opt& Time & Int Hyp  & Eigen Hyp & Opt & Time & Opt & Time   \\
\bottomrule

\multirow{5}{*}{MB} & 1 & 5 & 12.31 & 1.99 & 1 &  7.2 & 12.31 & 4.17 &  1.0  &  12.6 &12.31 & 0.13 & 12.31 &  0.53\\

& 2& 10 & 23.16 & 2.76 & 5.6 & 38.4 & 23.16 & 24.67  & 4.8  & 34.2 & 23.16 & 2.56 &  23.16 &  35.94\\

& 3 & 10 & 21.15 &  6.44 &  8.6  &  58.6 & 21.15 & 21.61 & 6.2 & 28.2 & 21.15 & 6.64 & 21.15 &  1.40\\

& 4 & 20 & 19.34 &  643.77 &  240.4   &   1326.4 & 19.34 &  300.46 &  26.2   &    112.4 & 19.34 & 21.47 & 19.34 &    6.80 \\

& 5 & 40 &  9.36  & 3600* & 258.4  & 543.6 & Inf &      1931.57 & 53.4   &   137 &  9.36 & 12.54 & 9.36 &  0.32\\

\hline

\multirow{5}{*}{MI} & 1 & 5 & 12.21 & 1.98   & 3 & 13.2 & 12.21 & 3.68 & 3 & 17 & 12.21 & 0.28 & 12.21  &  0.48\\

& 2& 10 & 23.16 &  1.59 & 6.4  & 27.6 & 23.16 &  9.83 & 6.4 &  27 & 23.16 &   2.11  & 23.16 &  1452.73\\

& 3 & 10 & 21.09  & 1.84 &  11.4  & 23 & 21.09 &   12.65 &  11.4 &   25.8 & 21.09 &  1.80 & 21.09  &   109.63\\

& 4& 20  &19.34 & 39.51 & 40.6  &   131.0 & 19.34 &  180.88 & 38.0 & 107.4 & 19.34 & 104.54 & 19.34 &      8.84\\

& 5 & 40 & 9.34 & 879.91 &  177.8    &  457.6 &  9.36& 1652.23  & 108.8 &    235.8 & 9.34 &  15.20 & 9.34 &  0.16\\

\bottomrule
\end{tabular}
}
\caption{Results for Problem \eqref{eq:lmpopt} when applying RPT-BB, RPT-SDP-BB, BARON, and SCIP for the mixed-binary (MB) and mixed-integer (MI) domains. \label{tab:numexp_lmp}}
\end{table} 

From Table \ref{tab:numexp_lmp}, we observe that in the mixed binary case, all methods are able to solve problem instances 1, 2, and 3 in less than a minute. Further, for problem instances 4 and 5 SCIP is the best performing method, followed by BARON, and then by RPT-BB. 

Moreover, in the mixed integer case, we observe a different behavior. For problem instances 2 and 3, both RPT-BB and BARON solve the problem in less than 3 seconds, while SCIP requires minutes to solve it. On the other hand, for problem instances 4 and 5, SCIP is the best performing method. Further, we observe that RPT-BB is faster than BARON on problem instance 4, while BARON is faster on problem instance 5.

\subsection{Dike height optimization} \label{sec:dike}
The next problem that we consider is the dike heightening problem from \cite{dike}, where the authors propose an optimization model to determine the optimal dike heightening in the Netherlands to prevent flooding. In this paper, we consider a special time-truncated version of the model described in \cite{bertsimas2023rptbb}:

\small
\begin{align} \label{eq:dho} 
            \begin{array}{cll}
            \displaystyle	\min_{\bm x, \bm h \in \mathbb{Z}_{+}^N } \underbrace{\sum_{k \in \mathcal{K}_0} \left(C + bx_k\right) \exp\left(\lambda h_k - \delta t_k\right)}_{\rm Investment \ costs} + \underbrace{ \sum\limits_{k \in \mathcal{K}_0} \frac{S_0}{\beta_{\delta}}\left(\exp\left(\beta_{\delta} t_{k+1}\right) - \exp\left(\beta_{\delta} t_k \right)\right)\exp\left(-\theta h_k \right)}_{\rm Expected \ damage \ costs} + \underbrace{\vphantom{\sum_{k \in K} f_k} \frac{S_0}{\delta}\exp\left(\beta_{\delta}T - \theta h_K \right)}_{\rm Future \ damage \ costs}.
            \end{array}
\end{align}
\normalsize

\noindent where $t$ represents the vector of time instances when the dike height is increased, $x_k$ is the integer height increment at time $t_k$, $h_k = \sum_{i=0}^k x_k$ is the total increment since $t_0$, and $C, b, \lambda, \delta, S_0, \beta_{\delta}, T,\text{and } \theta$, are constants. 

In \cite{bertsimas2023rptbb}, the authors consider a case of the dike height problem where the height increases are in the continuous domain. In this paper, we consider a simplification of that problem, where all height changes are fixed to be discrete in the integer domain.

We evaluate the performance of our algorithm on Problem \eqref{eq:dho}, using the instances from \cite{bertsimas2023rptbb}. We consider three different frequencies of dike height updates. In the first one, $t_{25}$,  we update the height of the dike every 25 years, in the second one, $t_{50}$, every $50$ years, and in the third one, $t_{ir}$, in irregular time increments. For each case, we examine a dike height problem with 10, 15, and 16 dike rings. We record the results in Table \ref{tab:numexp_dho}.

\begin{table}[H]
\begin{center}
\resizebox{\columnwidth}{!}{
 \begin{tabular}{cccccccccccccccccc}

\toprule 
 \multirow{2}{*}{} & \multirow{2}{*}{\#}& \multirow{2}{*}{$n_x$} & \multicolumn{3}{c}{\textbf{RPT-BB}}& \multicolumn{3}{c}{\textbf{RPT-SDP-BB}} &\multicolumn{2}{c}{\textbf{BARON}}&\multicolumn{2}{c}{\textbf{SCIP}}
\\ \cmidrule(r){4-6}\cmidrule(r){7-9}\cmidrule(r){10-11}\cmidrule(r){12-13}
&& & Opt& Time &Int Hyp  &Opt& Time &Int Hyp & Opt & Time & Opt & Time \\
\bottomrule
\multirow{3}{*}{$t_{25}$}  

&10& 12  & 61.31& 	$3600^*$	&14922&61.31 &	\textbf{5.77} &2 & 96.21 &	$3600^*$ &	61.31 &	8.04\\

& 15& 12 & 610.06 & $3600^*$	& 48919& 609.94 &	$\bm{3600^*}$	&4037&618.25 &	$3600^*$ &	609.94 &	$\bm{3600^*}$ \\

& 16 & 12 & 1269.68 & $3600^*$ &	46859& 1269.64	& \textbf{136.43} &	144 &2014125.46	&$3600^*$&1269.65	&$3600^*$  \\
\hline

\multirow{3}{*}{$t_{50}$}
& 10 & 6 & 155.5	& $3600^*$ &	13076 & 155.5 &	4.86 &2& 89.1 &	$3600^*$ &	155.5	& \textbf{2.35} \\

& 15 & 6  &545.26 &	$3600^*$ &	22839 & 545.26 &	5.97	&7& 545.66 & $3600^*$ &	545.26 & \textbf{3.23} \\

& 16 &6  &1100.15&	88.65	&4482 & 1100.15	&7.20&	10&	1116.28	& $3600^*$	&1100.15	& \textbf{2.51} \\
\hline
\multirow{3}{*}{$t_{ir}$}
& 10 & 5 &61.98 & $3600^*$	& 20897& 61.98 &	$3600^*$ &	7541& 96.85	& $3600^*$ &	61.98	& \textbf{3.1} \\

& 15 & 20 & 608.76	&$3600^*$ &	57410 & 608.76 &  \textbf{124.91}	&219& 624.53 &	$3600^*$	&608.77 &	$3600^*$\\

& 16 & 10 & 1268.18 & $3600^*$	& 58763 &1268.18 &	\textbf{110.38} &	229& 	2058.60
&$3600^*$&	1268.28	&$3600^*$ \\

\hline

\end{tabular}}
\end{center}
\caption{Results for Problem \eqref{eq:dho} when applying RPT-BB, RPT-SDP-BB, BARON, and SCIP, for dike rings 10, 15, and 16.\label{tab:numexp_dho}}%
\end{table}

We observe that while RPT-BB is only able to solve the $t_{50}$ case for dike ring $16$ to provable global optimality, when we add the SDP constraint to the RPT relaxation, the RPT-SDP-BB algorithm is able to solve 7 out of 9 problems to optimality within two minutes. Additionally, looking at the state-of-the-art solvers, we first see that BARON is unable to solve any problem to provable global optimality. However, we observe that SCIP solves some problems to provable optimality, doing so within a few seconds. 

Comparing the two best-performing approaches, RPT-SDP-BB and SCIP, we first see that both algorithms solve all instances in the $t_{50}$ case to global optimality within a few seconds, with SCIP being the fastest. Additionally, in the $t_{25}$ case we observe that RPT-SDP-BB outperforms SCIP for both dike ring 10 and dike ring 16. Notably, for dike ring 16 we see that SCIP is unable to solve the problem within the time limit of an hour. Finally, in the $t_{ir}$ case, we observe that SCIP outperforms RPT-SDP-BB for dike ring 10, while for dike rings 15 and 16, RPT-SDP-BB is significantly faster, with SCIP being unable to solve to provable optimality.  

\subsection{Comparison of branching strategies} 
\label{sec:branch_strat} 
In this section, we experiment with three different branching strategies to solve Problems \eqref{eq:smlt} and \eqref{eq:QCQPnconvex} for the mixed-binary and mixed-integer domains, and demonstrate the value for each one. The three branching strategies are as follows:
\begin{enumerate}
    \item First perform integer branching. Once all discrete variables are fixed, branch on the continuous variables using eigenvector branching. More specifically, we set the fractionality treshold $\delta$ equal to $0$, see also Step 3 in Section \ref{sec:method}.
    \item Perform integer branching only if there exists an integer variable whose value is more than $\delta = 0.1$ away from the closest integer. 
    \item Randomly select between integer and eigenvector branching with a fixed probability. More specifically, let $p$ be the probability of performing integer branching and $1-p$ the probability of using eigenvector branching. Then, if any of the values of the integer variables is not integer, we will perform an integer branch with probability $p$. For this strategy we experiment with $p = 0.8$.
\end{enumerate}
The results are given in Table \ref{tab:numexp_somol_strategies} and Table \ref{tab:numexp_qcqp_strategies}.





\begin{table}[H]
\resizebox{\columnwidth}{!}{

 \begin{tabular}{cccccccccccccccccc}
\toprule 
\multirow{2}{*}{} & \multirow{2}{*}{$\set X$} & \multirow{2}{*}{\#} & \multirow{2}{*}{$n_x$} &\multicolumn{4}{c}{\textbf{Original}}&\multicolumn{4}{c}{\textbf{Threshold 0.1}}
&\multicolumn{2}{c}{\textbf{Random 80/20}}
\\ \cmidrule(r){5-8}\cmidrule(r){9-12}\cmidrule(r){13-14}\cmidrule(r){15-16}
&&&&Opt& Time & Int Hyp  & Eigen Hyp&Opt& Time & Int Hyp  & Eigen Hyp & Opt & Time & Opt & Time \\
\bottomrule


\multirow{10}{*}{MB}& \multirow{5}{*}{$\set X_1$} & 1 &  $5$  & 10.93 &	8.80	&0&	0&  10.9266 & 7.41 &   0  &  0 & 10.93 &  7.16 &   0  &  0\\

&& 2 & $5$ & 89.52	& 0.02	&0&	0 	&89.52 & 0.02 &   0  &  0 & 89.52  & 0.02 & 0  &   0\\

&& 3 & $20$ & 244.51	& 73.30	&381 &	12	&	244.51 &  71.67 &   381  &      12 & 244.51 & 154.36 &    452 &       104\\

&& 4 &$10$ & 68.37 & 0.01 &	0 &	0	&	 68.37 & 0.01 & 0   &  0 &  68.37 &  0.01 &  0  & 0\\

&& 5  &$20$ & 2649.33 &	40.21  &	92	& 31 &	 2649.33 & 24.48  &  46   &  26 & 2649.33  & 13.72 &   22  &  11\\

\cline{2-16}

&\multirow{5}{*}{$\set X_2$}& 1 & $5$ & 10.93	&8.30 &	0 &	0	&	10.93  & 8.43 &   0   &  0 &  10.93 & 8.07   & 0   &  0 \\

&& 2 & $5$ & 89.52	& 0.03 &	0	&0		&89.52 & 0.03 &  0   &   0 & 89.52 & 0.03 &  0  & 0 \\

&& 3 &	20 &244.51	& 288.71 &	377&	12	&	244.51 &  301.73 &    377   &     12 & 244.51 & 539.09 &    382  &      87 \\

&& 4 &	10 & 31.18	&0.63&	7	&0		&31.18 & 0.75 &   6 &  2 & 31.18 &  0.61&   4     &     2\\

&& 5 &	20 & 1483.24&	155.84&	71	&47	&	1483.24  & 152.925  &  66     &    48	& 1483.24  & 279.657  &  118  &      79\\
\bottomrule


\multirow{10}{*}{MI} & \multirow{5}{*}{$\set X_1$}& 1 &  $5$ & 22.85	&8.19&	1&	0	&22.85 & 8.02 &    1   & 0 & 22.85 & 7.51  & 2   & 1\\
&& 2 & $5$ & 224.93	& 0.04  &	1&	0	&224.93 & 0.04 & 1       &   0 & 224.93 & 0.27 & 5   &  4\\
&& 3 & $20$ & 1069.34	& 47.78	&218	&0	&1069.34 & 47.67  &  218   &  0 &  1069.34 &  124.25 &  402  &   91\\
&& 4 &$10$ & 113.32&	0.09	&5&	0		& 113.32 &  3600*  &  0   & 6887 & 113.32  & 0.25 &  9  &        4\\
&& 5 &$20$ & 2987.78	& 5.17 &12	&0	&	2987.78  & 5.24 &  12     &    0 & 2987.78  & 5.10  & 12      &  0\\

\cline{2-16}
&\multirow{5}{*}{$\set X_2$}  & 1& $5$  & 13.69	&9.07&	3&	0	&	113.69 & 9.08  & 3       &   0 & 13.69 &  8.97 &  4      &    1 \\

&& 2& $5$ & 134.48	&0.23 &	2&	0	&	119.40 & 3600* &  1     &     10787 & 134.48 & 0.20  & 2     &     0\\

&& 3 &	20 &809.78&	94.31	&125	&0	&	809.78  & 96.78 & 125     &   0 &   809.78  & 274.16 &  236        & 50\\

&& 4 &	10 & 31.91	&1.05&	11	&0	&	31.91 & 1.06 &  11     &    0 &  31.91 & 2.26 &  19      &   5 \\

&& 5 &	20 & 1593.48&	3600*	&2406	&0	&	1593.48 & 3600* & 2385    &   8 & 1593.48 & 3600*  & 1737      & 393\\
\bottomrule

\end{tabular}
}
\caption{Results of using RPT-BB with three different branching strategies on the mixed-binary (MB) and mixed-integer (MI) domains for Problem \eqref{eq:smlt}. \label{tab:numexp_somol_strategies}}
\end{table} 

From Table \ref{tab:numexp_somol_strategies}, we observe that all three strategies perform similarly on the mixed-binary domain with set $\set X_1$. However, for set $\set X_2$, the strategies using thresholds of 0 and 0.1 clearly outperform the random strategy. In the mixed-integer domain, the strategy with threshold 0 shows superior performance overall. Specifically, for $\set X_1$, instance 3, the threshold-based strategies (0 and 0.1) outperform the random baseline. For instance 4, however, the strategies with threshold 0 and the random approach significantly outperform the threshold 0.1 strategy, which fails to prove optimality within the time limit. These results indicate that the threshold 0 strategy yields the most consistent performance across instances in this domain.


\begin{table}[h]
\resizebox{\columnwidth}{!}{
 \begin{tabular}{ccccccccccccccccc}
\toprule 
 \multirow{2}{*}{$\set X$} & \multirow{2}{*}{\#}&  \multirow{2}{*}{$n_x$} & \multicolumn{4}{c}{\textbf{Threshold 0}}&\multicolumn{4}{c}{\textbf{Threshold 0.1}}
&\multicolumn{4}{c}{\textbf{Random 80/20}}
\\ \cmidrule(r){4-7}\cmidrule(r){8-11}\cmidrule(r){12-15}
&&& Opt& Time & Int Hyp  & Eigen Hyp&Opt& Time & Int Hyp  & Eigen Hyp & Opt & Time & Int Hyp & Eigen Hyp   \\
\bottomrule

\multirow{10}{*}{MB} & 1 & 20 & 247.76  &	7.68  &	2 &	3	&	247.76  &  7.66   & 1    &  2&		247.76  &  7.80  &  1    &   4\\

& 2& 20 & 737.76 &	0.85 &	2	& 3	&737.76  &  0.79 &  1      &    2   &737.76  &  0.90 &  2   &   3\\

& 3 & 10 & 2706.69	& 0.06	& 0	& 0	&	2706.69 &   0.07 & 0     &     0	 & 2706.69  &  0.06 & 0   &     0\\

& 4 & 50 &  77333.33	&3.93 &	3&	0	&	 77333.33  &  7.86  &  1      &    4 & 77333.33  &  7.70   & 1     &     4 \\

& 5 & 100 & 572687.27	&6.40 &	0&	0	&	572687.27 &6.47   &  0    &    0	& 572687.27 & 6.49  &  0           & 0 \\

& 6 & 8 & 408.00	& 0.30	& 0	& 0	&	408.00    &  0.30 &  0      &    0  &	408.00    &  0.30  & 0      &    0 \\

& 7& 12 &  282.06&	0.65	&9&	0		&282.06 &   0.65 &  8      &    1 & 282.06  &  1.66  &  17 &        2\\

& 8 & 16 & 521.35&	0.21 &	6&	0		&521.35 & 3600* & 6     &     3002 &	521.35  & 0.22 &  6      &    0\\

& 9 & 30 & 2149.32&	4.98	&38&	0	& 2149.32 & 22.44 &  39    &     6 &  2149.32  &  12.22  &   31   &   6\\

& 10 & 40 & 2421.69	&174.83	&506&	0		&2421.69  &182.22  & 504    &    6  & 2421.69  &  860.70  &  1160       & 279\\
\hline 
\multirow{10}{*}{MI} & 1 & 20 & 375.84&	7.87 & 4&	0	&	375.841   &  8.14  &  4   &    0	& 375.84 &    10.61  &  11     &    7\\

& 2& 20 & 865.84	& 0.79	&4	&0	&865.84  &  0.80 &  4     &     0   & 865.84  &  0.81  &  4      &    0\\

& 3 & 10 & 6888.78	&0.06&	0	&0	&	6888.78  &  0.06 &  0       &   0	& 6888.78  &  0.06 & 0   &       0\\

& 4 & 50 & 98331.20 & 2.43 &	2 &	0	& 98331.20  &  2.25 &  2    &  0  &	98331.20   & 2.16  & 2   &   0\\

& 5 & 100 & 774403.58	&  90.64 &	2 &	0	& 774403.58 & 3600*   & 1  &  30 & 774403.58 & 279.72  &  3 &       1\\

& 6 & 8 & 1717.80 &	0.41 &	0 &	0	&	1717.80  &   0.31  & 0   &  0 & 1717.80  &  0.31  & 0  &        0\\

& 7& 12 & 4505.24	& 0.07 &	2 &	0	&	4505.24    & 0.08 & 2  &  0  & 4505.24  &  0.07  & 2          & 0\\

& 8 & 16 & 17063.66 &	0.16 &	3 &	0	&	17063.66  &  0.26  &  2     &     3  & 17063.66   & 0.16  & 3          & 0\\

& 9 & 30 & 42821.70	&$3600^*$	&9029	&3151	&	44851.60   & 3600*  &  8017    &   2834 & 44900.10  &  3600*     & 3083    &   1095\\

& 10 & 40 & 23914.80	&$3600^*$&	1482&	1185	&	23914.80 &   3600*  &  1346   &  1194 & 25497.40  &  3600* &    1245    &   1262\\
\hline
\end{tabular}
}
\caption{Results of using RPT-BB with three different branching strategies on the mixed-binary (MB) and mixed-integer domains for Problem \eqref{eq:QCQPnconvex}.}\label{tab:numexp_qcqp_strategies}%
\end{table}

Table \ref{tab:numexp_qcqp_strategies} further confirms this observation. The strategy with threshold 0 consistently outperforms the alternatives in scenarios requiring substantial branching.  Most notably, in the mixed-binary domain, for instances $8, 9$, and $10$, both other strategies perform at least as bad as the strategy with $\delta =  0$, with the strategy with threshold 0.1 not able to prove optimality within the time limit for problem instance 8. 
Additionally, for the mixed-integer domain, for problem instance 5 the strategy with threshold 0.1 is not able to prove optimality within the time limit of one hour. Moreover, for problem instances 9 and 10, all strategies are not able to prove optimality within the time limit, with the strategy with threshold 0 finding the best upper bound.


In summary, the strategy that performs binary branching followed by eigenvector-based branching emerges as the most effective among those we evaluated.

\section{Discussion and conclusion} \label{sec:conclusions}
We develop a method for globally solving mixed integer nonconvex optimization problems involving SLC functions. We leverage the RPT framework along with the continuous relaxation of the integer variables, to obtain a convex relaxation of the original nonconvex problem, while introducing additional variables and constraints. We then develop a branch and bound scheme, leveraging a combination of eigenvector branching with integer branching rules, to obtain the global optimal solution of the initial problem. In the numerical experiments, we demonstrate that the proposed method stands well against the current state of the art global optimization methods for mixed integer nonconvex optimization problems. Our approach outperforms BARON for many problem instances in terms of both computational time and optimal objective value. We observe that SCIP continues to lead as the state-of-the-art solver. In contrast, BARON shows limitations in handling increased nonconvexity or integer complexity, often failing to scale efficiently. Our approach demonstrates strong performance on challenging instances, successfully solving problems to global optimality that SCIP and BARON are unable to solve within the one-hour time limit.

An alternative to our current approach is to reverse the order in which RPT and variable relaxation is applied. Instead of first relaxing the integer variables and then applying RPT to the resulting nonconvex continuous problem to obtain a convex relaxation, we may instead apply RPT directly to the original mixed-integer optimization problem. This yields a convex mixed-integer formulation that can be handled by a solver such as SCIP. Given that our results show that SCIP continues to lead as the state-of-the-art solver, particularly for problems where variables are predominantly binary or general integer, we can leverage its strengths for branching on the integer variables. Once the integer structure has been fully resolved through branching, we then can apply eigenvector-based branching on the continuous variables (basically setting $\delta = 0$). An important question is therefore which approach, applying RPT after relaxing the integer variables or immediately applying RPT and using SCIP to solve the resulting convex mixed-integer formulation, yields better overall performance. Moreover, in the integer and mixed-integer setting, one key advantage of our approach is that the added branching constraints can be directly incorporated into the RPT formulation to tighten the convex relaxation further, an effect that is also reflected in our numerical results. Investigating this question represents a promising direction for future research and could provide valuable insights into selecting the most effective strategy across different variable domains.

\bibliography{Masterbib}
\bibliographystyle{plainnat}

\newpage

\end{document}